\documentclass[12pt]{article}

\usepackage{amsmath}
\usepackage{amsfonts}
\usepackage{amsthm}
\usepackage{hyperref}
\usepackage{enumitem}
\usepackage{array}
\usepackage[nottoc,numbib]{tocbibind}

\newcommand{\N}{\mbox{$\mathbb{N}$}}
\newcommand{\nN}{\mbox{$(n\in\N)$}}

\begin{document}

%1111111111111111111111111111111111111111111111111111111111111111111111111111

\author{Attila Losonczi}
\title{On infinite versions of the prisoner problem}

\date{\today}%{8 September 2021}

\newtheorem{thm}{\qquad Theorem}[section]
\newtheorem{prp}[thm]{\qquad Proposition}
\newtheorem{lem}[thm]{\qquad Lemma}
\newtheorem{cor}[thm]{\qquad Corollary}
\newtheorem{rem}[thm]{\qquad Remark}
\newtheorem{ex}[thm]{\qquad Example}
\newtheorem{df}[thm]{\qquad Definition}
\newtheorem{prb}{\qquad Problem}

\newtheorem*{thm2}{\qquad Theorem}
\newtheorem*{cor2}{\qquad Corollary}
\newtheorem*{prp2}{\qquad Proposition}

\maketitle

\begin{abstract}

\noindent

We investigate some versions of the famous 100 prisoner problem for the infinite case, where there are infinitely many prisoners and infinitely many boxes with labels. In this case, many questions can be asked about the admissible steps of the prisoners, the constraints they have to follow and also about the releasing conditions. We will present and analyze many versions and cases. In the infinite case, the solutions and methods require mainly analysis rather than combinatorics.

\noindent
\footnotetext{\noindent
AMS (2020) Subject Classifications:  40A99, 05A05\\

}

\end{abstract}

%\tableofcontents

%-------------------------------------------------------------------------------------------------------------------------------------------Introduction----------------------------------
\section{Introduction}

First, let us recall the original ''100 prisoner'' problem.

There are 100 prisoners and in a separate room, there are 100 boxes. The guard has labels from 1 to 100, and he puts the labels into the boxes randomly in a way that each box contains exactly one label. Then, one by one, the prisoners enter the room and each of them can open 50 boxes. After one prisoner finishes, the boxes are closed and no information is allowed to be sent back. If all prisoners find his label, then all of them are released, otherwise none. The question is what is the best strategy? I.e. what strategy can maximize the probability of releasing?

Many authors have investigated the 100 prisoner problem and also its several variations and got surprising answers, see e.g. \cite{gm}, \cite{fs},  \cite{ab},  \cite{adi}, \cite{wc},  \cite{w},  \cite{gs}, \cite{ckp}, \cite{lsw}. 

\smallskip

The aim of this paper is to investigate the infinite version of this problem, namely, when there are infinitely many prisoners and infinitely many boxes and infinitely many labels. Here, many natural questions can be asked. In this paper, we start to investigate these kinds of problems. 
\par Here, in the introduction, we just sketch the basics of our problems, for more precise and detailed description see section \ref{s1} and subsequent sections. Roughly, our basic setting is the following.
\vspace{-0.3pc}\begin{itemize}\setlength\itemsep{-0.2em}
\item Infinity is countable infinity, i.e. the prisoners, boxes, labels are numbered by the natural numbers.
\item We suppose that the guard \textit{somehow} puts the labels into the boxes such that each box contains one label and all label is used.
\item Each box has a price that must be paid to open it. The prisoners get the same amount as the sum of the prices of the boxes. Their task is to divide it among them.
\end{itemize}

From here, there may be many directions that can be followed. In this paper, we will investigate the case when each prisoner does exactly the same than in the solution of the original problem, i.e. opens the box with the same label he has, then opens the box which has the label that he just found, etc. (a so-called pointer-following algorithm). He does that till his money lasts. Furthermore, we assume that the guard arranged the labels so that all cycles are finite. For releasing conditions we will investigate two cases. All prisoners are released if
\vspace{-0.3pc}\begin{itemize}\setlength\itemsep{-0.2em}
\item infinitely many of them find his label
\item all prisoners except finitely many of them find his label.
\end{itemize}

It is important to note that unlike the original problem, we do not investigate randomness, i.e. it is not a probability problem in its current form.

\smallskip

In this paper, we investigate five types of problems which can be grouped into two main sets.

%------------------------------------------------------------------------------------------------------------------------------Versions---------------------------
\section{Versions of the problem}

%------------------------------------------------------------------------------------------------------------------------------Version 1.a---------------------------
\subsection{Version 1.a}\label{s1}

First we are going to investigate the following problem.

\vspace{-0.6pc}\begin{itemize}\setlength\itemsep{-0.3em}
\item There are infinitely many prisoners, labeled by natural numbers.
\item There are infinitely many boxes, labeled by natural numbers as well.
\item The prison guard puts the labels of the prisoners into the boxes in a way that each box contains exactly one number and every number is in a box, i.e. it is a permutation of $\N$. The prisoners do not know this permutation.
\item Furthermore, all cycles in the permutation are finite, i.e. if we open an arbitrary box, look at the number inside and open the box labeled by that number, and follow that procedure, then in finitely many steps we will and up with the number of the originally opened box.
\item Each box has a weight or price that has to be paid in order to open it. The sum of all prices equals to $1$. The prisoners know the price of each box.
\item The group of prisoners gets $1$ as a full amount. They have to divide this amount among them. Before the first prisoner enters the room, they decide how much each of them gets.
\item The prisoners go into the room of boxes one by one.
\item A prisoner has to pay for opening each box he wants to open. A prisoner can open as many boxes as he wants provided he can cover the costs of the opened boxes.
\item If a prisoner is in the room where the boxes are, neither he can send back any information, nor he can send back the not used amount.
\item All prisoner does exactly the same thing: First open the box with the same label that the prisoner has. Then look at the number inside, and open the box labeled by that number. And so on.  (pointer-following algorithm) He does that until his money lasts and can cover opening the next box.
\item When a prisoner finishes and leaves the room, all boxes get closed.
\item Finally, all prisoners are released if infinitely many of them found his label. Otherwise none of them.
\item Question: Under which circumstances do the prisoners have a strategy for getting released for any possible permutation of the labels, and what is that strategy?
\end{itemize}

Note that the prisoners does exactly the same steps than in the classic solution of the original 100 prisoner problem. As the cycles are finite, each prisoner has a chance to find his label. In this problem, \textbf{their strategy is not what they do in the room of the boxes, instead, the strategy is ''only'' how they divide the original amount among them}.

Hence the question is the following. Knowing the prices of the boxes, can the prisoners divide the full price in a way that infinitely many of them find his label using the above mentioned pointer-following algorithm for any admissible permutation of the labels. In other words, how does the box price distribution determines a strategy, if there is any?

Very important to emphasize that for a given box price distribution, we are looking for a strategy that works for ALL admissible permutations of the labels.

\smallskip

We note again that this is not a probability problem (unlike the original 100 prisoner problem), we do not investigate the probability of getting released.

\smallskip

Later, when we have a necessary and sufficient condition for the existence of a strategy, we will see that many of the above defining conditions can be weakened (see remarks \ref{r2}, \ref{r3}, \ref{r4}, \ref{r5}).

\subsubsection*{Notions and notations used in the sequel}

We will denote the set of positive integers by $\N$.

\smallskip

If $\sigma:\N\to\N$ is a permutation, a cycle will be denoted by $(n_1,n_2,\dots,n_k)$ where $n_i\in\N\ (1\leq i\leq k)$, $n_i\ne n_j\ (1\leq i,j\leq k, i\ne j)$ and $\sigma(n_i)=n_{i+1}\ (1\leq i\leq k-1)$ and $\sigma(n_k)=n_1$.

\smallskip

We will denote the permutation set by the guard by $\sigma:\N\to\N$, that is simply a bijection of $\N$. We call $\sigma$ a \textbf{finite-cycle-permutation} if all cycles are finite.

\smallskip

We will denote the price of the $n^{th}$ box by $p_n$, and the amount that the $n^{th}$ prisoner gets by $a_n$, i.e. the strategy of the prisoners is the sequence $(a_n)$. We require that $p_n\geq 0,\ a_n\geq 0$ and $\sum_{n=1}^{\infty}p_n=\sum_{n=1}^{\infty}a_n=1$.

\smallskip

The \textbf{price of a cycle} $(n_1,\dots,n_k)$ is $\sum_{i=1}^kp_{n_i}$.

%------------------------------------------------------------------------------------------------------------------------------Solutions---------------------------
\subsubsection{Solutions}

First, we add a trivial necessary condition for $(a_n)$ being a strategy.

\begin{prp}If $(a_n)$ is a strategy, then for infinitely many $n\in\N$, $a_n\geq$ the price of the cycle in which $n$ is. Especially $a_n\geq p_n$ has to hold for infinitely many $n\in\N$. 
\end{prp}
\begin{proof}In order to find his label, the $n^{th}$ prisoner has to open all boxes belonging to his cycle.
\end{proof}

Let us investigate a basic price distribution first. Suppose that the $n^{th}$ box has price $p_n=\frac{1}{2^n}$. 
\par Then a strategy can be the following. The first prisoner gets $a_1=0$, and the $n^{th}$ prisoner ($n>1$) gets amount $a_n=\frac{1}{2^{n-1}}$. Clearly $\sum_{n=1}^{\infty}p_n=\sum_{n=1}^{\infty}a_n=1$.
\par Let us verify that it works. We claim that the prisoner with label equals to the least member of each cycle will find his label (except one cycle with least element $1$). As the length of the cycles are finite, there are infinitely many cycles, therefore infinitely many prisoner will find his label, hence all will be released.
\par Let $\sigma:\N\to\N$ be the finite-cycle-permutation set by the guard. Let $\{n_1,\dots,n_k\}\subset\N$ be a cycle, i.e. $\sigma(n_j)=n_{j+1}$ for $1\leq j\leq k-1$, and $\sigma(n_k)=p(n_1)$. We can suppose that $n_1<n_j$ if $2\leq j\leq k$, i.e. $n_1$ is the least member of the cycle. We also suppose that $1<n_1$. We show that prisoner $n_1$ will find label $n_1$. Prisoner $n_1$ has amount $a_{n_1}=\frac{1}{2^{n_1-1}}$. But that amount is enough to open all boxes with label $\geq n_1$, because $\sum_{i=n_1}^{\infty}p_i=\sum_{i=n_1}^{\infty}\frac{1}{2^i}=\frac{1}{2^{n_1-1}}$. And all members of the cycle is among those boxes, hence he will able to open all of them, and when he opens box $n_k$, then he will find his own label which is $n_1$.

We can simply generalize this result by noticing that with suitable choice of amounts, the rest of the prisoners can open a ''tail'' of all boxes.

\begin{thm}\label{t1}If $\sum_{n=1}^{\infty}np_n<\infty$, then there is a strategy.
\end{thm}
\begin{proof}Let $\tilde p_n=\sum_{m=n}^{\infty}p_m$. If $\sum_{n=1}^{\infty}np_n<\infty$, then we get that $\sum_{n=1}^{\infty}np_n=\sum_{n=1}^{\infty}\tilde p_n<\infty$. Choose $m\in\N, m>1$ such that $\tilde p=\sum_{n=m}^{\infty}\tilde p_n<1$.
\par Let the strategy be the following. Let
\[a_n=\begin{cases}
1-\tilde p&\text{if }n=1\\
0&\text{if }1<n<m\\
\tilde p_n&\text{if }n\geq m.
\end{cases}\]
We claim that a prisoner with label equals to the least member of a cycle with having least member at least $m$ will find his label.
\par Let $\sigma:\N\to\N$ be the finite-cycle-permutation set by the guard. Let $\{n_1,\dots,n_k\}\subset\N$ be a cycle. Suppose that $n_1<n_j$ if $2\leq j\leq k$. We also suppose that $m\leq n_1$. We show that prisoner $n_1$ will find label $n_1$. Prisoner $n_1$ has amount $a_{n_1}=\tilde p_{n_1}=\sum_{m=n_1}^{\infty}p_m$. But clearly that amount is enough to open all boxes with label $\geq n_1$. And all members of the cycle is among those boxes, hence he will able to open all of them, and when he opens box $n_k$, then he will find his own label which is $n_1$.
\end{proof}

\begin{rem}\label{r1}From the proof, we can see that theorem \ref{t1} remain valid if the prisoners get any (small) finite amount (not 1 as originally). 
\end{rem}

We can go on generalization by noticing that the order of the boxes may not count.

\begin{thm}\label{t2}If there is a permutation $\delta:\N\to\N$ such that ${\sum_{n=1}^{\infty}np_{\delta(n)}<\infty}$, then there is a strategy.
\end{thm}
\begin{proof}The prisoners know the prices of the boxes and if $\sum_{n=1}^{\infty}np_n=\infty$, but there is a permutation $\delta:\N\to\N$ such that $\sum_{n=1}^{\infty}np_{\delta(n)}<\infty$, then they can rearrange the order of the boxes in their heads according to $\delta$. And then they apply the strategy in theorem \ref{t1} for $(p_{\delta(n)})$.
\end{proof}

A question arises immediately if how this property of a price distribution can be changed by a permutation.

\begin{prp}Let $(p_n)$ be a sequence such that $p_n\geq 0\ \nN$ and suppose that there are infinitely many indexes $n_k$ such that $p_{n_k}>0$. Then there is a permutation $\delta:\N\to\N$ such that $\sum_{n=1}^{\infty}np_{\delta(n)}=\infty$.
\end{prp}
\begin{proof}Let $(n_k)$ be the sequence of indexes for which $p_{n_k}>0$. Then find a strictly increasing sequence $(j_k)$ of natural numbers such that ${j_kp_{n_{2k-1}}>1}$. Obviously it can be found. Then simply move the element $p_{n_{2k-1}}$ to the position $j_k$, and move the rest of the sequence $(p_n)$ such that they fill the remaining places.
\end{proof}

Therefore there are sequences for which a permutation can change that kind of sum. We show that there are also sequences for which that sum does not change by a permutation.

\begin{thm}\label{t4}Let $p_n=\frac{1}{n^2}\ \nN$ and $\delta:\N\to\N$ be a permutation. Then $\sum_{n=1}^{\infty}np_{\delta(n)}=\infty$.
\end{thm}
\begin{proof}First note that it is enough to prove that $\sum_{n=1}^{\infty}\delta(n)p_n=\infty$ for all permutation $\delta$, because switching to $\delta^{-1}$ we get that 
\[\sum_{n=1}^{\infty}np_{\delta(n)}=\sum_{n=1}^{\infty}\delta^{-1}(n)p_n.\]
\par Let $\delta:\N\to\N$ be a permutation. We will show more, namely, that for all $m\in\N$ we have
\[\sum_{n=1}^{m}\frac{1}{n}=\sum_{n=1}^{m}np_n\leq\sum_{n=1}^{\infty}\delta(n)p_n,\]
i.e. the original ordering (or identity permutation) gives the minimum.
\par Let $m\in\N$ be fixed. For $n\in\N, n\leq m$, let $\alpha_n=\delta(n)-n$ the size of the step at $n$, which can be either positive, negative or zero. Clearly
\[\delta(n)p_n=\delta(n)\frac{1}{n^2}=(n+\alpha_n)\frac{1}{n^2}=\frac{1}{n}+\alpha_n\frac{1}{n^2}.\]
Let
\[t_n=\begin{cases}
\frac{1}{n}&\text{if }\alpha_n=0\\
\frac{1}{n}+\sum\limits_{i=1}^{\alpha_n}\frac{1}{(n+i)^2}&\text{if }\alpha_n>0\\
\frac{1}{n}+\sum\limits_{i=0}^{-\alpha_n+1}-\frac{1}{(n-i)^2}&\text{if }\alpha_n<0.
\end{cases}\]
It is easy to see that in each cases
\[t_n\leq\frac{1}{n}+\alpha_n\frac{1}{n^2}=\delta(n)p_n\]
holds, because when $\alpha_n\ne 0$, then we substitute each $\frac{1}{n^2}$ by a smaller term (there are $\alpha_n$ copies of $\frac{1}{n^2}$, and we have the same number of new terms). Hence if we substitute $\delta(n)p_n$ by $t_n$, then we make the sum smaller. 
\par Now we are going to define the so-called semi-cycles. There are three types of them. 
\vspace{-0.4pc}\begin{enumerate}\setlength\itemsep{-0.3em}
\item We call $(n)$ a semi-cycle if $n\leq m,\ \delta(n)=n$.
\item $(n_1,n_2,\dots,n_k)$ is a semi-cycle if the following holds.
\vspace{-0.5pc}\begin{itemize}\setlength\itemsep{-0.2em}
\item $n_l\leq m\ (1\leq l\leq k)$
\item $\delta(n_l)=n_{l+1}\ (1\leq l\leq k-1)$
\item $\delta^{-1}(n_1)>m$
\item $\delta(n_k)>m$.
\end{itemize}
\vspace{-0.5pc}In words, we came outside of $\{1,\dots,m\}$, and we step till we are out of $\{1,\dots,m\}$ again.
\item $(n_1,n_2,\dots,n_k)$ is a semi-cycle if the following holds. 
\vspace{-0.5pc}\begin{itemize}\setlength\itemsep{-0.2em}
\item $n_l\leq m\ (1\leq l\leq k)$
\item $\delta(n_l)=n_{l+1}\ (1\leq l\leq k-1)$
\item $\delta(n_k)=n_1$
\item $n_l<n_1\ (1\leq l\leq k-1)$.
\end{itemize}
\vspace{-0.5pc}In words, it is a full cycle inside $\{1,\dots,m\}$ such that $n_1$ is the greatest term. (Actually this case contains the first one.)
\end{enumerate}
\par It can be readily seen $\{1,\dots,m\}$ can be decomposed into disjoint union of semi-cycles.
\par In a semi-cycle, let us make the substitution mentioned above, i.e. substitute all $\delta(n)p_n$ by $t_n$. It is easy to see that all negative terms $-\frac{1}{i^2}$ will be canceled by a positive one inside the semi-cycle. Therefore we get that
\[\sum_{i=1}^{k}\delta(n_i)p_{n_i}=\sum_{i=1}^{k}n_{i+1}p_{n_i}\geq\sum_{i=1}^{k}t_{n_i}\geq\sum_{i=1}^{k}\frac{1}{n_i},\]
i.e. the original sum of the terms in the semi-cycle is greater than or equal to the reciprocals of the terms in the semi-cycle. If we do it for all semi-cycles, then we get the claim.
\end{proof}

Now we can show the counterpart of theorem \ref{t2}.

\begin{thm}\label{t3}Let $(p_n)$ be a sequence such that $p_n\geq 0\ \nN$ and $\sum_{n=1}^{\infty}p_n=1$, and for every permutation $\delta:\N\to\N$, $\sum_{n=1}^{\infty}np_{\delta(n)}=\infty$ holds. Then there is no strategy.
\end{thm}
\begin{proof}First let us analyze what ''\textit{there is no strategy}'' means. For a strategy $(a_n)$, we require that for ALL finite-cycle-permutations of the labels, infinitely many prisoners find his label. A strategy depends on the price distribution only. Hence ''\textit{there is no strategy}'' is equivalent with that for every $(a_n)$ there is a finite-cycle-permutation $\sigma:\N\to\N$ for which only finitely many prisoners will find his label. In other words, if the prisoners produce $(a_n)$, then for that, the guard can make a finite-cycle-permutation of the labels such that only finitely many prisoners will find his label.
\par Let $(a_n)$ be the strategy of the prisoners. 
\par If there are zeros among the elements of $(a_n)$, then distribute amount 1 among them somehow (does not matter how). I.e. modify the values of such $a_n$-s by adding altogether 1 to them. E.g. if there are infinitely many zeros, then modify the zeros to $\frac{1}{2},\frac{1}{4},\frac{1}{8},\dots$. To simplify the notation, let us denote this new sequence by $(a_n)$ as well. Clearly the sum changes to $\sum_{n=1}^{\infty}a_n=2$.
\par Let us reorder $(a_n)$ to become descending, i.e. find a permutation ${\delta:\N\to\N}$ such that $(a_{\delta(n)})$ is descending. To make to notation simpler, we can assume that $(a_n)$ is descending, i.e. we solve the problem for descending $(a_n)$ with the corresponding $(p_n)$, and after solving it, apply $\delta^{-1}$ for both sequences to get back the original ones.
\par Let
\[s_n=\sum_{i=n}^{\infty}p_i\ \ \nN.\]
We call $n\in\N$ good, if $s_n>a_n$, otherwise we call $n$ bad.
\par We claim that there are infinitely many natural numbers that are good. Suppose the contrary, i.e. there is an $N\in\N$ such that $n\geq N$ implies that $n$ is bad. It means that $s_n\leq a_n\ (n\geq N)$. Therefore
\begin{equation}\label{eq1}\sum_{n=N}^{\infty} s_n\leq\sum_{n=N}^{\infty}a_n\leq 2\end{equation}
but
\[\sum_{n=N}^{\infty} s_n=\sum_{n=N}^{\infty}\sum_{i=n}^{\infty}p_i=\sum_{m=1}^{\infty}mp_{N-1+m}=\sum_{m=N}^{\infty}mp_m-\sum_{m=N}^{\infty}(N-1)p_m\]
\[=\sum_{m=1}^{\infty}mp_m-\sum_{m=1}^{N-1}mp_m-(N-1)\sum_{m=N}^{\infty}p_m=\infty\]
because $\sum_{m=1}^{\infty}mp_m=\infty$ and $\sum_{m=N}^{\infty}p_m\leq 1$. We get that $\sum_{n=N}^{\infty} s_n$ is infinite and $\leq 1$ at the same time which is a contradiction.
\par Now we are ready to formalize how the guard can make a suitable finite-cycle-permutation of the labels.
\vspace{-0.6pc}\begin{enumerate}\setlength\itemsep{-0.3em}
\item First he changes the zeros in $(a_n)$ by adding altogether 1 to them (as we did above).
\item Then he reorders $(a_n)$ to become descending (as we did above).
\item He finds the first good $n$.
\item He finds consecutive elements $n,n+1,\dots,n+k$ such that 
\vspace{-0.6pc}\begin{enumerate}\setlength\itemsep{-0.3em}
\item $\sum\limits_{l=0}^{k}p_{n+l}>a_n$
\item $n+k+1$ is good.
\end{enumerate}
Then he closes his cycle, i.e. this is his cycle: $(n,n+1,\dots,n+k)$.
\item He repeats step (4) for $n+k+1$ instead of $n$.
\end{enumerate}
Obviously (4a) will happen for some $k\geq 0$ because $n$ is good. And at the same time, (4b) can also be fulfilled as there are infinitely many good numbers.
\par Clearly in each such cycle, the first member will not find his label by (4a) (if the original $a_n$ was 0, then it is trivial). But as $(a_n)$ is descending, none of the rest members in the cycle will (again, if some member was originally 0, then it is trivial). Hence after the first such cycle, no one will find his label.
\end{proof}

\begin{rem}\label{r2}From the proof, we can see that theorem \ref{t3} remains valid if the prisoners get any finite amount (not 1 as originally), because we just used it in inequality \ref{eq1} and its consequence.
\par Together with remark \ref{r1}, we can say that the answer (if there is a strategy) does not depend on how much (finite) amount the prisoners get altogether.\qed
\end{rem}

We can now merge theorems \ref{t2} and \ref{t3} into one.

\begin{cor}\label{c1}There is a strategy if and only if there is a permutation $\delta:\N\to\N$ such that $\sum_{n=1}^{\infty}np_{\delta(n)}<\infty$.\qed
\end{cor}

\begin{rem}\label{r3}We can weaken the condition regarding how the guard is supposed to create his permutation. It can be readily seen that corollary \ref{c1} remains valid if we just require that the guard creates permutation such that it contains infinitely many finite cycles (and we allow that it also contains infinite cycles as well).\qed
\end{rem}

\begin{rem}\label{r4}Moreover we can drop the ''finite-cycle-permutation'' condition completely, if we also change how the box opening happens. When a prisoner enters the room, he shows his money to the guard and he points to a box, and he asks the guard to open all boxes belonging to the full cycle which contains that box, provided his money covers all boxes. If his money does not cover all, then he has to say which boxes he wants to open one-by-one. Again, in this case, corollary \ref{c1} remains valid.\qed
\end{rem}

\begin{rem}\label{r5}As a consequence of remark \ref{r2}, we can say that corollary \ref{c1} remains valid if we modify one more condition, namely, we allow the prisoners to send back the not used amount. Moreover everyone can send back all his amount. Because it acts like we just double the original full amount (which is 1).\qed
\end{rem}

We can characterize those sequences for which this peculiar sum in corollary \ref{c1} is infinite for all permutations. First, we characterize it for positive sequences.

\begin{thm}\label{t5}Let  $(p_n)$ be a sequence such that $p_n>0\ \nN$. Then $\sum_{n=1}^{\infty}np_{\delta(n)}=\infty$ holds for every permutation $\delta:\N\to\N$ if and only if $\sum_{n=1}^{\infty}np_{\delta(n)}=\infty$ holds when $\delta$ is a permutation which orders the sequence to be descending, i.e. $p_{\delta(n+1)}\leq p_{\delta(n)}\ \nN$.
\end{thm}
\begin{proof}\par The necessity is obvious. We prove the sufficiency.
\par Suppose that $\delta$ is a permutation which orders the sequence to be descending. Hence (by reordering the sequence) we can assume that $p_{n+1}\leq p_n$ and $\sum_{n=1}^{\infty}np_{n}=\infty$.
\par Then we have to modify the definition of $t_n$ in the proof of theorem \ref{t4} only, and all the rest of the proof will work here as well. 
\par Recall from \ref{t4} that $\alpha_n=\delta(n)-n$, and $\delta(n)p_n=(n+\alpha_n)p_n=np_n+\alpha_np_n$.
So in our case let
\[t_n=\begin{cases}
np_n&\text{if }\alpha_n=0\\
np_n+\sum\limits_{i=1}^{\alpha_n}p_{n+i}&\text{if }\alpha_n>0\\
np_n+\sum\limits_{i=0}^{-\alpha_n+1}-p_{n-i}&\text{if }\alpha_n<0.
\end{cases}\]
Then all rest of the proof of \ref{t4} can be repeated thereafter.
\end{proof}

Before the general case, we need a little preparation.

\begin{df}Let $(a_n)$ be a non-negative sequence. We call $(a_n)$ \textbf{quasi-descending}, if $0<a_n,\ n<m$ implies that $a_m\leq a_n$ holds.\qed
\end{df}

\begin{prp}Let $(a_n)$ be a non-negative sequence. Then $(a_n)$ can be reordered to become quasi-descending.
\end{prp}
\begin{proof}Let us order the strictly positive elements in a descending order, then add the remaining 0 elements somehow.
\end{proof}

\begin{thm}\label{t9}Let  $(p_n)$ be a sequence. Let $(q_n)$ be the sequence that we get when we omit the zero terms from $(p_n)$. Then the following statements are equivalent.
\vspace{-0.6pc}\begin{enumerate}\setlength\itemsep{-0.3em}
\item $\sum_{n=1}^{\infty}np_{\delta(n)}=\infty$ holds for every permutation $\delta:\N\to\N$.
\item $\sum_{n=1}^{\infty}np_{\delta(n)}=\infty$ holds when $\delta$ is a permutation which orders the sequence to be quasi-descending.
\item ${\sum_{n=1}^{\infty}nq_{\delta(n)}=\infty}$ holds for every permutation $\delta:\N\to\N$.
\item ${\sum_{n=1}^{\infty}nq_{\delta(n)}=\infty}$ holds when $\delta$ is a permutation which orders the sequence $(q_n)$ to be descending.
\end{enumerate}
\end{thm}
\begin{proof}If $p_n>0$, then let $\alpha(n)$ denote the index of the number $p_n$ in the sequence $(q_n)$, in other words let $q_{\alpha(n)}=p_n$ for $p_n>0$. Clearly $\alpha:\N\to\N$ is injective and descending.

\smallskip

The equivalence of 3 and 4 is in theorem \ref{t5}. Clearly the implications $1\Rightarrow 2$ and $3\Rightarrow 4$ are obvious.

\smallskip

$3\Rightarrow 1$: If there is a permutation $\delta$ on $(p_n)$, then it naturally determines a permutation $\sigma$ on the positive elements of $(p_n)$ i.e. on $(q_n)$. But the sum gets smaller: $\infty=\sum_{n=1}^{\infty}nq_{\sigma(n)}\leq\sum_{n=1}^{\infty}np_{\delta(n)}$.

\smallskip

$4\Rightarrow 2$: The same argument works with noting that a quasi-descending ordering on $(p_n)$ is a descending ordering on $(q_n)$. 

\smallskip

$1\Rightarrow 3$: Let $\delta:\N\to\N$ be a permutation (applied for $(q_n)$). We define a permutation $\sigma$ of $(p_n)$. Let $\sigma:\N\to\N$ be the following. If $p_n>0$, then let $\sigma(n)=2\delta(\alpha(n))\ \nN$, i.e. we do exactly the same like $\delta$ does, but on the set of even numbers only. For odd numbers, we somehow arrange the zeros of $(p_n)$, that is how we define $\sigma$. Clearly we get that $\infty=\sum_{n=1}^{\infty}np_{\sigma(n)}=2\sum_{n=1}^{\infty}nq_{\delta(n)}$. Hence $\sum_{n=1}^{\infty}nq_{\delta(n)}=\infty$.

\smallskip

$2\Rightarrow 4$: The previous argument works with noting that $\sigma$ is a quasi-descending ordering of $(p_n)$ if $\delta$ is a descending ordering of $(q_n)$. 
\end{proof}

Now we can simply derive if we omit the zeros in $(p_n)$, then we get an equivalent sequence regarding strategy.

\begin{cor}Let  $(p_n)$ be given. Let us denote by $(q_n)$ the sequence that we get when we omit the zeros from $(p_n)$. Then there is strategy for $(p_n)$ if and only if there is one for $(q_n)$.
\end{cor}
\begin{proof}Corollary \ref{c1} and theorem \ref{t9}.
\end{proof}

Analyzing the proofs of \ref{t4} and \ref{t5}, we can realize that some parts can be applied to sequences with $\sum_{n=1}^{\infty}np_{\delta(n)}<\infty$ for a permutation $\delta:\N\to\N$.

\begin{prp}Let  $(p_n)$ be a sequence such that $p_n>0\ \nN$ and $\sum_{n=1}^{\infty}np_{\delta(n)}<\infty$ for a permutation $\delta:\N\to\N$. Let $\sigma:\N\to\N$ be a permutation such that $p_{\sigma(n+1)}\leq p_{\sigma(n)}$. Then for any permutation $\delta:\N\to\N$ and any $m\in\N$ we have that
\[\sum_{n=1}^{m}np_{\sigma(n)}\leq\sum_{n=1}^{m}np_{\delta(n)}.\]
Especially $\sum_{n=1}^{\infty}np_{\sigma(n)}<\infty$.
\end{prp}
\begin{proof}Apply the corresponding parts of theorems \ref{t4} and \ref{t5}.
\end{proof}

We can partially apply the proof of theorem \ref{t3} for the case when $\sum_{n=1}^{\infty}np_{\delta(n)}<\infty$ for a given permutation $\delta:\N\to\N$.

\begin{prp}Let $\delta:\N\to\N$ be a permutation such that $a_{\delta(n+1)}\leq a_{\delta(n)}$. If $\sum_{n=1}^{\infty}np_{\delta(n)}=\infty$, then the guard can make a finite-cycle-permutation such that the prisoners will not be released.
\end{prp}
\begin{proof}Apply the corresponding parts of theorems \ref{t3}.
\end{proof}

\begin{cor}If there is a strategy $(a_n)$, then for any permutation $\delta:\N\to\N$ such that $a_{\delta(n+1)}\leq a_{\delta(n)}$, $\sum_{n=1}^{\infty}np_{\delta(n)}<\infty$ has to hold.\qed
\end{cor}

Now we investigate the case when the lengths of the cycles are bounded.

\begin{thm}\label{t6}Suppose that the set of lengths of all cycles is bounded and this bound is known by the prisoners. Then there is a strategy.
\end{thm}
\begin{proof}Let the bound be $k\geq 1$.
Then there is $m\in\N$ such that ${k\sum_{i=m}^{\infty}p_i<1}$. The strategy is $(a_n)$ where $a_n=0$ if $n<m$, and $a_n=kp_n$ if $n\geq m$. If we take a cycle which does not contain element from $\{1,\dots,m\}$, then in the cycle, prisoner $n$ can find its label for which $p_n$ is the greatest in the cycle, because his amount $a_n$ is greater than or equal to the price of the cycle.
\end{proof}

%------------------------------------------------------------------------------------------------------------------------------Version 1.b---------------------------
\subsection{Version 1.b}\label{s1b}

Everything is the same as in Version 1.a (see section \ref{s1}), except that we change the releasing condition. Now all prisoners are released if each except finitely many of them finds his label. This is a significant strengthening of the condition as the next statements will show.

\smallskip

In this section, when we say ''strategy'' then it means that it refers to strategy according to Version 1.b.

\begin{prp}In Version 1.b, there is no strategy in general.
\end{prp}
\begin{proof}As usual, let the prices of the boxes be $(p_n)$ and let the strategy of the prisoners be $(a_n)$.
\par Let the first cycle be $(1,2,\dots,k_1)$ where $k_1=\left\lceil\frac{1}{p_1}\right\rceil+1$. Then there must be at least one $j\in\{1,\dots,k_1\}$ such that $a_j\leq\frac{1}{k_1}$ because $\sum_{n=1}^{\infty}a_n=1$. But then \mbox{prisoner $j$} will not be able to open box 1 because $a_j<p_1$, hence he will not find his label.
\par Suppose that we have already ordered the numbers $1,\dots,m-1$ into cycles. Then our next cycle will be $(m,m+1,\dots,m+k_m)$ where $k_m=\left\lceil\frac{1}{p_{m}}\right\rceil$. Again there must be at least one $j\in\{m,m+1,\dots,m+k_m\}$ such that $a_j\leq\frac{1}{k_m+1}$ because $\sum_{n=1}^{\infty}a_n=1$. But again \mbox{prisoner $j$} will not be able to open box $m$ because $a_j<p_m$, hence he will not find his label.
\par We can conclude that in each cycle, there is at least one prisoner who does not find his label.
\end{proof}

However under not too strict conditions, there is strategy. 
\par First note that as our condition in Version 1.b is stronger than the condition in Version 1.a, in order to have a strategy, it is necessary that there is a permutation $\delta:\N\to\N$ such that $\sum_{n=1}^{\infty}np_{\delta(n)}<\infty$ (see theorem \ref{t3}).

\begin{thm}Let $\sum_{n=1}^{\infty}np_{\delta(n)}<\infty$ for a suitable permutation ${\delta:\N\to\N}$, and suppose that the diameter of cycles is bounded in $(\delta(n))$ (not in $(n)$) and this bound is known by the prisoners (together with $\delta$). Then there is strategy.
\end{thm}
\begin{proof}The prisoners' first action is that they reorder the boxes in order to get $\sum_{n=1}^{\infty}np_{n}<\infty$ (to simplify the notation, we omit $\delta$).
\par Let the diameter of cycles $\leq d\in\N$. As $\sum_{n=1}^{\infty}np_{n}<\infty$, there is an $m\in\N$ such that
\[\sum_{n=m}^{\infty}s_n<1\text{ where }s_n=\sum_{i=n}^{\infty}p_i.\]
Then the strategy is the following. The first $m+d$ prisoner gets $0$ as amount. Then the prisoner numbered $m+d+k\ (k\in\N)$ gets amount $s_{m+k}$.
\par Let us verify if it works. For prisoner numbered $m+d+k$, when he opens box numbered $m+d+k$, he enters a cycle with diameter less than or equal to $d$. Hence the minimum member in the cycle can be numbered as $m+k$. But the sum of prices of all later boxes equals to $s_{m+k}$. And as there are finitely many boxes in the cycle, he can open all of them and finds his label.
\end{proof}

Here, knowing the bound of the length of all cycles does not help.

\begin{thm}\label{t7}Suppose that the set of lengths of all cycles is bounded and this bound, $k$ is known by the prisoners. Then there is no a strategy even for $k=2$.
\end{thm}
\begin{proof}Let $(a_n)$ be the amount distribution of the prisoners. Then the guard does the following. For $n=1$, he finds $n_1\in\N$ such that $a_{n_1}<p_1$, and then he create cycle $(1,n_1)$. If he finished with $m$ cycles, then he takes the minimum element $l\in\N$ such that $l$ is not in the previously created cycles. Then he finds $n_l\in\N$ such that $n_l$ is not in the previously created cycles and $a_{n_l}<p_l$. Then he create cycle $(l,n_l)$. Obviously it is good, because we got infinitely many $2$-length long cycles, and at least one member will not find his label.
\end{proof}

%------------------------------------------------------------------------------------------------------------------------------Version 1.c---------------------------
\subsection{Version 1.c}

It is a modification of Version 1.b (see section \ref{s1b}) with the following two differences. 
\vspace{-0.6pc}\begin{itemize}\setlength\itemsep{-0.3em}
\item The guard does not close the boxes after a prisoner leaves the room. I.e. the next prisoner entering the room gets extra information as all boxes opened by the previous prisoners will be revealed to him immediately.
\item The prisoners can decide the order how they enter the room, i.e. they do not have to follow their labels as the order for entering the room. (Until now, this condition was not important, now, it is.)
\end{itemize}

The releasing condition remains the same as in Version 1.b. 

In this case, knowing the bound of the length of all cycles does help.

\begin{thm}\label{t8}Suppose that the set of lengths of all cycles is bounded and this bound, $k$ is known by the prisoners. Then there is a strategy.
\end{thm}
\begin{proof}First, the prisoners order $(p_n)$ to be descending, i.e. $p_n\leq p_{n+1}$ (omitting $\delta$ to simplify the notation). Then we can copy the method described in theorem \ref{t6}. I.e.  there is $m\in\N$ such that ${k\sum_{i=m}^{\infty}p_i<1}$. The strategy is $(a_n)$ where $a_n=0$ if $n<m$, and $a_n=kp_n$ if $n\geq m$.
Then the smallest member of each cycle with $n>m$ will be able to open all boxes in the cycle, and as the boxes will remain open, all later member in the cycle will find his label.
\end{proof}

%------------------------------------------------------------------------------------------------------------------------------Version 1.c---------------------------
\subsection{Version 1.d}

It is a modification of Version 1.b (see section \ref{s1b}) with the difference that the guard tells the prisoners all cycles, more precisely all members of each cycle (but not the order inside of course).

\begin{thm}There is no strategy in general.
\end{thm}
\begin{proof}Let $(1,\dots, k_1)$ be the first cycle created by the guard such that there is $i\in\{1,\dots,k_1\}$ such that $p_i>0$ and $k_1p_i>1$. If the guard has created cycles with members less than or equal to $m\in\N$, then his next cycle is $(m+1,\dots,m+k_m)$ such that there is $i\in\{m+1,\dots,m+k_m\}$ such that $p_i>0$ and $k_mp_i>1$. Evidently this can be done.
\par As the prisoners gets altogether 1 as amount, there must be at least one prisoner in each cycle who cannot open all boxes of his cycle.
\end{proof}

\begin{thm}Suppose that the set of lengths of all cycles is bounded and this bound, $k$ is known by the prisoners. Then there is a strategy.
\end{thm}
\begin{proof}There is $m\in\N$ such that ${k\sum_{i=m}^{\infty}p_i<1}$. If a prisoner's label is $n\geq m$, and all members of his full cycle $\geq m$, then he gets the amount which is the price of his cycle. Everyone else gets 0. Clearly this strategy works, and the full used amount is less than or equal to ${k\sum_{i=m}^{\infty}p_i<1}$.
\end{proof}

%------------------------------------------------------------------------------------------------------------------------------Version 2.a---------------------------
\subsection{Version 2.a}

There are two differences to Version 1.a which are:
\vspace{-0.6pc}\begin{enumerate}\setlength\itemsep{-0.3em}
\item The prices of the boxes are fixed, moreover the sum of all prices is infinite. The price of the $n^{th}$ box is $\frac{1}{n}$.
\item A strategy of the prisoners is a sequence of non-negative numbers $(a_n)$ (there is no restriction for the sum of $(a_n)$).
\end{enumerate}
Hence the releasing condition is unchanged, i.e. all prisoners are released if infinitely many of them find their labels.

Now we are going to analyze some strategies of the prisoners.

%----------------------------------------------------------
\subsubsection*{Strategy: $a_n=1$ for all $n\in\N$}

It does not work always. Just create cycles such that the price of each cycle is greater than 1.

%----------------------------------------------------------
\subsubsection*{Strategy: $a_n=\sum\limits_{i=1}^n\frac{1}{i}\ \nN$}

Clearly the last member of each cycle will find his label, hence all prisoners will be released.

%----------------------------------------------------------
\subsubsection*{Strategy: There is a fixed $k\in\N$. If $n<k$, then $a_n=0$. If $n\geq k$, then $a_n=\sum\limits_{i=k}^n\frac{1}{i}$}

There are only finitely many cycles which contains element that is smaller than or equal to $k$. If a cycle is not of that type, then the last member of the cycle will find his label, hence all prisoners will be released.

%----------------------------------------------------------
\subsubsection*{Strategy: There is a fixed $K\in\mathbb{R}$. Let $a_n=\log n - K$}

This case can be simply reduced to the previous one. For $K$, there is a $k\in\N$ such that $\sum_{i=k}^n\frac{1}{i}\leq\log n-K$. Because $\sum_{i=1}^n\frac{1}{i}\leq\log n-1$, and then $K+1\leq\sum_{i=1}^{k+1}\frac{1}{i}$ is enough.

%----------------------------------------------------------
\subsubsection*{Strategy: $a_n=c\sum\limits_{i=1}^n\frac{1}{i}\ \nN$ for fixed $c<1$}

We claim that it does not work always. To show that it is enough to prove the following.
\begin{prp}For all $k\in\N$ there is $n\geq k$ such that $c\sum\limits_{i=1}^n\frac{1}{i}<\sum\limits_{i=k}^n\frac{1}{i}$.
\end{prp}
\begin{proof}Clearly
\[\sum\limits_{i=1}^n\frac{1}{i}<\frac{1}{c}\sum\limits_{i=k}^n\frac{1}{i}\]
\[\sum\limits_{i=1}^{k-1}\frac{1}{i}<\left(\frac{1}{c}-1\right)\sum\limits_{i=k}^n\frac{1}{i}.\]
As $\left(\frac{1}{c}-1\right)>0$, a suitable $n$ can be found.
\end{proof}
Now if we want to put $k$ into a new cycle, then find $n$ according to the proposition, and let the cycle be $(k, k+1,\dots,n)$.

%------------------------------------------------------------------------------------------------------------------------------Version 2.b---------------------------
\subsection{Version 2.b}

We just modify the releasing condition to that all prisoners are released if each except finitely many of them finds his label.

Then there is no strategy. For fixed $n$, just create a cycle ${(n,n+1,\dots,n+k)}$ such that $\sum_{i=n}^{n+k}\frac{1}{i}>a_n$. 

%---------------------------------------------------------------------------------------------------------------------------------------------------------the bibliography -----------------

\medskip

{\footnotesize

%------------------------------------------------------------------------------------------------------------------------------------------------------------address---------------
\noindent

\noindent email: alosonczi1@gmail.com\\
}
\end{document}